\documentclass[11pt]{amsart}
\usepackage{cras}

\usepackage[english,francais]{babel}
\usepackage{amsfonts, amsmath,amssymb}
\usepackage{hyperref}
\usepackage[all]{xy}
\usepackage[lite]{amsrefs}
\usepackage{graphicx}
\usepackage{wrapfig}
\usepackage{pstricks, pstricks-add}

\newcommand{\ringO}{\mathcal{O}}

\newcommand{\C}{{\mathbb{C}}}

\newcommand{\Z}{{\mathbb{Z}}}

\newcommand{\rationals}{{\mathbb{Q}}}

\newcommand*{\Homol}{\operatorname{H}}

\newcommand{\Afour}{\mathcal{A}_4}
\newcommand{\Sthree}{\mathcal{S}_3}

\newcommand{\circlegraph}{ 
\begin{pspicture}       (-0.3,-0.3)(0.3,0.3)
                        \pscircle(0,0.0){0.25}
                        \psdots(0.25,0.0)
\end{pspicture} }

\newcommand{\edgegraph}{ 
\begin{pspicture}(-0.3,-0.3)(0.3,0.3)
        \psdots(-0.2,0.0)
        \psdots(0.2,0.0)
        \psline(-0.2,0.0)(0.2,0.0)
\end{pspicture} }

\newtheorem{theorem}{Theorem}[section]
\newtheorem{lemma}[theorem]{Lemma}
\newtheorem{e-proposition}[theorem]{Proposition}

\newtheorem{e-definition}[theorem]{Definition\rm}

\newtheorem{theoreme}{Th\'eor\`eme}[section]

\def\og{\leavevmode\raise.3ex\hbox{$\scriptscriptstyle\langle\!\langle$~}}
\def\fg{\leavevmode\raise.3ex\hbox{~$\!\scriptscriptstyle\,\rangle\!\rangle$}}

\begin{document}

\selectlanguage{english}
\title{Homology and $K$-theory of the \mbox{Bianchi} groups}

\selectlanguage{english}
\author{Alexander D. Rahm}
\email{Alexander.Rahm@Weizmann.ac.il}
\urladdr{http://www.wisdom.weizmann.ac.il/\char126rahm/}

\address{Department of Mathematics, Weizmann Institute of Science, Rehovot, Israel}

\issueinfo{+++}{+++}{January}{2011}
\date{\today}
\seriesinfo{Th\'eorie des groupes}{Group Theory}

\begin{abstract}
\selectlanguage{english}
We reveal a correspondence between the homological torsion of the \mbox{Bianchi} groups and new geometric invariants, which are effectively computable thanks to their action on hyperbolic space.
We use it to explicitly compute their integral group homology and equivariant $K$-homology. 
By the Baum/Connes conjecture, which holds for the \mbox{Bianchi} groups, we obtain the $K$-theory of their reduced $C^*$-algebras in terms of isomorphic images of the computed $K$-homology.

We further find an application to Chen/Ruan orbifold cohomology.
\end{abstract}

\selectlanguage{francais}
\ftitle{Homologie et $K$-th\'eorie des groupes de \mbox{Bianchi}}
\begin{fabstract} 
Nous mettons en \'evidence une correspondance entre la torsion homologique des groupes de \mbox{Bianchi} et de nouveaux invariants g\'eom\'etriques, calculables gr\^ace \`a leur action sur l'espace hyperbolique. Nous l'utilisons pour  calculer explicitement leur homologie de groupe \`a coefficients entiers et leur $K$-homologie \'equivariante. 
En cons\'equence de la conjecture de Baum/Connes, qui est v\'erifi\'ee  pour ces groupes, nous obtenons la $K$-th\'eorie de leurs C*-alg\`ebres r\'eduites en termes d'images isomorphes de la $K$-homologie calcul\'ee.
Nous trouvons d'ailleurs une application \`a la cohomologie d'orbi-espace de Chen/Ruan.

\end{fabstract}
\maketitle

\selectlanguage{francais}
\section{Version fran\c{c}aise abr\'eg\'ee}

Nous \'etudions la g\'eom\'etrie d'une certaine classe de groupes arithm\'etiques (les groupes de Bianchi), \`a travers une action propre sur un espace contractile.
Nous acc\'edons \`a leur homologie de groupe et leur $K$-homologie \'equivariante.
En plus de d\'etail, consid\'erons un corps de nombres quadratique imaginaire $\rationals(\sqrt{-m})$, o\`u $m$ est un entier positif ne contenant pas de carr\'e. Soit~$\ringO_{-m}$ son anneau d'entiers.
Les \emph{groupes de \mbox{Bianchi}} sont les groupes $\mathrm{SL_2}(\ringO_{-m})$. Les groupes de \mbox{Bianchi} peuvent \^etre consid\'er\'es cruciaux pour l'\'etude d'une classe plus large de groupes, les groupes \emph{Kleiniens}, qui ont d\'ej\`a \'et\'e \'etudi\'es par Henri \mbox{Poincar\'e}~\cite {Poincare}.
En fait, chaque groupe Kleinien arithm\'etique non-cocompact est commensurable avec un groupe de \mbox{Bianchi} \cite{MaclachlanReid}. Un \'eventail d'informations sur les groupes de \mbox{Bianchi} peut \^etre trouv\'e dans les monographies \cites{Fine, ElstrodtGrunewaldMennicke, MaclachlanReid}.
Ces groupes agissent d'une mani\`ere naturelle sur l'espace hyperbolique \`a trois dimensions, qui est isomorphe \`a l'espace sym\'etrique qui leur est associ\'e.
Le noyau de cette action est le centre $\{ \pm 1 \}$ des groupes, ce qui rend utile l'\'etude du quotient par le centre, $\mathrm{PSL_2}(\ringO_{-m})$, que nous appellerons \'egalement un groupe de Bianchi.
En 1892, Luigi \mbox{Bianchi}~\cite{Bianchi} a calcul\'e des domaines fondamentaux pour cette action pour quelques uns de ces groupes. Un tel domaine fondamental est de la forme d'un poly\`edre hyperbolique (\`a quelques sommets manquants pr\`es), et nous l'appellerons le \emph{poly\`edre fondamental de \mbox{Bianchi}}.
Le calcul du poly\`edre fondamental de \mbox{Bianchi} a \'et\'e impl\'ement\'e sur ordinateur pour tous les groupes de Bianchi~\cite{BianchiGP}.
Les images sous $\mathrm{SL_2}(\ringO_{-m})$ des faces de ce poly\`edre munissent l'espace hyperbolique d'une structure cellulaire.
Pour mieux observer la g\'eom\'etrie locale, nous passons au \emph{complexe cellulaire raffin\'e}, que nous obtenons en subdivisant cette structure cellulaire jusqu'\`a ce que les stabilisateurs dans $\mathrm{SL_2}(\ringO_{-m})$ fixent les cellules point par point.
Nous allons exploiter ce complexe cellulaire de diff\'erentes mani\`eres, afin de cerner des aspects diff\'erents de la g\'eom\'etrie de ces groupes.

\subsection*{Homologie de groupes} 

Un invariant essentiel des groupes est leur homologie (d\'efinie par exemple dans~\cite{Brown}). Nous pouvons la calculer pour les groupes de Bianchi en nous servant du complexe cellulaire raffin\'e. \`A ce fin, nous utilisons la suite spectrale \'equivariante de \mbox{Leray}/\mbox{Serre} qui part de l'homologie des stabilisateurs d'un ensemble repr\'esenta\-tif de cellules, et qui converge vers l'homologie du groupe de \mbox{Bianchi}. Nous  pr\'ecisons dans la proposition~\ref{non-Euclidean results}, l'homologie enti\`ere de  $\mathrm{PSL_2}(\ringO_{-m})$ dans les cas $m=$ 19, 43, 67 et 163, qui constituent tous les cas d'anneaux principaux non-Euclidiens.
Les cas d'anneaux Euclidiens
 sont d\'ej\`a connus de~\cite{SchwermerVogtmann}. Des r\'esultats r\'ecents pour des cas de groupe de classe d'id\'eaux non-trivial se trouvent dans~\cite{RahmFuchs}.
Nous remarquons que dans les quatre cas d'anneaux principaux non-Euclidiens, la torsion dans l'homologie de $\mathrm{PSL_2}(\ringO_{-m})$, est du m\^eme type d'isomorphisme. 
Pour comprendre ceci, nous consid\'erons, pour un nombre premier~$\ell$, le sous-complexe de l'espace d'orbites, compos\'e des cellules ayant des stabilisateurs qui contiennent des \'el\'ements d'ordre~$\ell$.
Nous l'appellerons le \emph{sous-complexe de $\ell$--torsion}.
L'\'enonc\'e suivant traite la mani\`ere dans laquelle son type d'hom\'eomorphisme d\'etermine la suite spectrale \'equivariante. Il est d\'emontr\'e par la r\'eduction des sous-complexes de torsion effectu\'ee dans~\cite{RahmThesis}.
Cette technique utilise le lemme~\ref{geometricRigiditytheorem} pour cerner le type de stabilisateur d'un sommet $v$ qui est adjacent \`a exactement deux ar\^etes dont les stabilisateurs admettent de la $\ell$--torsion.
Ensuite, ces deux ar\^etes et $v$ sont remplac\'es par une seule ar\^ete. Le th\'eor\`eme~\ref{actionSurLesAxes} et quelques informations homologiques sur les groupes finis en question sont utilis\'es  pour v\'erifier que les morphismes induits en homologie produisent les m\^emes termes sur la deuxi\`eme page de la suite spectrale \'equivariante qu'avant le remplacement.

\begin{theoreme}[cf. theorem \ref{invarianceUnderReduction}] \label{invarianceSousReduction}
La partie $\ell$--primaire de la deuxi\`eme page de la suite spectrale \'equivariante convergeante vers l'homologie des groupes de \mbox{Bianchi} d\'epend, hors de sa  ligne inf\'erieure, seulement du type d'hom\'eomorphisme du sous-complexe de $\ell$--torsion.
\end{theoreme}

Dans tous les cas d'anneaux principaux non-Euclidiens, les sous-complexes de 2--torsion, respectivement de 3--torsion, sont hom\'eomorphes, ce qui explique les r\'esultats de la proposition~\ref{non-Euclidean results}.
Derri\`ere le th\'eor\`eme ci-dessus,  il y a la correspondance suivante entre les sous-groupes cycliques non-triviaux des stabilisateurs des sommets, et les lignes g\'eod\'esiques autour desquelles ils effectuent une rotation. Il convient d'appeler ces lignes des \emph{axes de rotation}.

\begin{theoreme}[cf. theorem \ref{actionOnAxes}] \label{actionSurLesAxes}
Soit $v$ un sommet quelconque dans l'espace hyperbolique. 
L'action de son stabilisateur sur l'ensemble des axes de rotation passant par $v$,
induite par l'action du groupe de \mbox{Bianchi}, 
est \'equivalente \`a l'action par conjugaison de ce stabilisateur sur ses sous-groupes cycliques non-triviaux.
\end{theoreme}

Une \'etude cas par cas \cite{RahmThesis} pour tous les six types de sous-groupes finis dans les groupes Bianchi nous permet de d\'eduire du th\'eor\`eme~\ref{actionSurLesAxes}, le lemme~\ref{geometricRigiditytheorem} utilis\'e pour obtenir le th\'eor\`eme~\ref{invarianceSousReduction}. Des exemples pour le th\'eor\`eme~\ref{invarianceSousReduction} sont donn\'es pour la 3--torsion homologique de trente-six groupes de Bianchi sur le tableau~\ref{3torsionsubgraphs}.

\subsection*{\emph{K}-th\'eorie}
Avec l'information sur l'action des groupes de \mbox{Bianchi}, que nous obtenons en nous servant des \'enonc\'es et m\'ethodes d\'ecrits ci-dessus, nous pouvons calculer  l'homologie de Bredon de groupes de \mbox{Bianchi}, et en d\'eduire leur $K$-homologie \'equivariante. Des r\'esultats sont pr\'esent\'es dans le th\'eor\`eme~\ref{K-homology}. 
En cons\'equence de la conjecture de Baum/Connes, qui est v\'erifi\'ee pour les groupes de \mbox{Bianchi} \cite{JulgKasparov}, nous obtenons la $K$-th\'eorie des $C^*$-alg\`ebres r\'eduites des groupes de \mbox{Bianchi} comme images isomorphes.


\selectlanguage{english}
\section{Introduction}

We study the geometry of a certain class of arithmetic groups (the Bianchi groups) by means of a proper action on a contractible space. This helps to determine their group homology and their equivariant \mbox{$K$-homology.}
In more detail, we denote by $\rationals(\sqrt{-m})$, with $m$ a square-free positive integer, an imaginary quadratic number field, and by $\ringO_{-m}$ its ring of integers. 
The \emph{\mbox{Bianchi} groups} are the groups $\mathrm{SL_2}(\ringO_{-m})$. 
The \mbox{Bianchi} groups may be considered as a key to the study of a larger class of groups, the \emph{Kleinian} groups, which dates back to work of Henri Poincar\'e~\cite{Poincare}.
In fact, each non-cocompact arithmetic Kleinian group is commensurable with some \mbox{Bianchi} group~\cite{MaclachlanReid}.  
A wealth of information on the \mbox{Bianchi} groups can be found in the monographs \cites{Fine, ElstrodtGrunewaldMennicke, MaclachlanReid}.
These groups act in a natural way on hyperbolic three-space, which is isomorphic to the symmetric space SL$_2(\C)/$SU$_2$ associated to them.
The kernel of this action is the centre $\{ \pm 1 \}$ of the groups.
Thus it is useful to study the quotient of a Bianchi group by its centre, namely $\mathrm{PSL_2}(\ringO_{-m})$, which we also call a Bianchi group.
In 1892, Luigi \mbox{Bianchi}~\cite{Bianchi} computed fundamental domains for this action when $m$ = 1, 2, 3, 5, 6, 7, 10, 11, 13, 15 and 19.
Such a fundamental domain has the shape of a hyperbolic polyhedron (up to a missing vertex at certain cusps, which represent the ideal classes of $\ringO_{-m}$), so we will call it the \emph{\mbox{Bianchi} fundamental polyhedron}. 
The computation of the \mbox{Bianchi} fundamental polyhedron has been implemented for all Bianchi groups~\cite{BianchiGP}.
The images under $\mathrm{SL_2}(\ringO_{-m})$ of the facets of this polyhedron equip hyperbolic three-space with a cell structure.
In order to view clearly the local geometry, we pass to the \emph{refined cell complex}, which we obtain by subdiving this cell structure until the cell stabilisers fix the cells pointwise. 
We will see how to exploit this cell complex in different ways, in order to see different aspects of the geometry of these groups.

\subsection*{Group homology}
An essential invariant of groups is their homology (defined for instance in~\cite{Brown}). We can compute it for the Bianchi groups using the refined cell complex and the equivariant Leray/Serre spectral sequence which starts from the homology of the stabilisers of representatives of the cells, 
and converges to the group homology of the \mbox{Bianchi} groups.
We will now state the results for simple integer coefficients in the cases $m = 19,$ $43,$ $67$ and $163$, which are the non-Euclidean principal ideal domain cases. In contrast to these, the Euclidean principal ideal domain cases are already known from~\cite{SchwermerVogtmann}. For some results in class number 2, see~\cite{RahmFuchs}.
The virtual cohomological dimension of the Bianchi groups is 2. In degrees strictly above 2, we express their homology in terms of the following Poincar\'e series at the primes $\ell = 2$ and $\ell = 3$: 
$$P^\ell_m(t) := \sum\limits_{q \thinspace = \thinspace 3}^{\infty} \dim_{\mathbb{F}_\ell} \Homol_q \left(\text{PSL}_2\bigl(\mathcal{O}_{-m}\bigr);\thinspace \Z/\ell \Z \right)\thinspace t^q.$$ 
These two primes are the only numbers which occur as orders of non-trivial finite elements of ${\rm PSL}_2(\ringO_{-m})$. So it has been shown \cite{RahmThesis} that the integral homology of these groups is, in all the mentioned degrees, a direct sum of copies of $\Z/2 \Z$ and $\Z/3 \Z$.
\begin{e-proposition} \label{non-Euclidean results}
The integral homology of ${\rm PSL}_2(\ringO_{-m})$, for \mbox{ $m \in \{ 19, 43, 67, 163\}$,}
is of isomorphism type $ \Homol_q({\rm PSL}_2(\ringO_{-m}); \thinspace \Z) \cong $
$ \begin{cases}
\Z^{\beta_1 -1} \oplus \Z/4\Z \oplus \Z/2\Z \oplus \Z/3\Z, & q = 2, \\
\Z^{\beta_1}, & q = 1,  \\ 
\end{cases} $ 
\\
where  \scriptsize
$\begin{array}{l|ccccc}
m        &  19 & 43 & 67 & 163 \\
\hline 
\beta_1  &  1   & 2  &  3 & 7   \\
\end{array}$ \normalsize
gives the Betti number $\beta_1$, and is in all higher degrees a direct sum of copies of $\Z/2 \Z$ and $\Z/3 \Z$, with the number of copies specified by the Poincar\'e series \begin{center}
$P^2_m(t) = \frac{-t^3(t^3 - 2t^2 + 2t - 3)}{(t-1)^2 (t^2 + t + 1 ) }$ and
$P^3_m(t) = \frac{-t^3(t^2 - t + 2)}{(t-1)(t^2+1)}$. \end{center}
\end{e-proposition}

We remark that in these four cases, the torsion in the integral homology of $ {\rm PSL}_2(\ringO_{-m})$ is of the same isomorphism type.
To understand this, we consider, for a prime $\ell$, the subcomplex of the orbit space consisting of the cells with elements of order $\ell$ in their stabiliser. We call it the \emph{$\ell$--torsion subcomplex}.
The following statement on how its homeomorphism type determines the equivariant spectral sequence is proven by the reduction of the torsion subcomplex carried out in~\cite{RahmThesis}. 
This technique uses lemma \ref{geometricRigiditytheorem} to determine the possible type of stabiliser of a vertex $v$ with exactly two adjacent edges which have $\ell$--torsion in their stabilisers.
Then these two edges, together with $v$, are replaced by a single edge; and theorem~\ref{actionOnAxes} as well as some homological information about the finite groups in question are used to check that the induced morphisms on homology produce the same terms on the second page of the equivariant spectral sequence as before the replacement.

\begin{theorem} \label{invarianceUnderReduction}
The $\ell$--primary part of the second page of the equivariant spectral sequence converging to the integral homology of ${\rm PSL}_2(\ringO_{-m})$ depends outside the bottom row only on the homeomorphism type of the $\ell$--torsion subcomplex.
\end{theorem}

\begin{figure} \caption{Results for the 3--torsion homology, expressed in $ P^\ell_m(t)$}

\label{3torsionsubgraphs} 

\scriptsize
\begin{tabular}{c|c|c}
$m$ specifying the \mbox{Bianchi} group  &  \begin{tabular}{c} {3-torsion subcomplex,} \\ {homeomorphism type} \end{tabular} 
& $ P^3_m(t)$ 
\\
\hline &  & \\
\begin{tabular}{c} 
{2, 5, 6, 10, 11, 15, 22, 29, 34, 35,} \\
{46, 51, 58, 87, 95, 115, 123, 155, } \\
{159, 187, 191, 235, 267} \end{tabular}
& \circlegraph &$\frac{-2t^3}{t-1}$ 
\\ &  &\\
7, 19, 43, 67, 139, 151, 163  & 
\edgegraph
 &$\frac{-t^3(t^2 - t + 2)}{(t-1)(t^2+1)}$ 
\\&  &\\
13, 37, 91, 403, 427 &  
\edgegraph \edgegraph  & 
$2  \left( \frac{-t^3(t^2-t+2)}{(t-1)(t^2+1)} \right)$ 
\\&  &\\
39 & \circlegraph \edgegraph & $\frac{-2t^3}{t-1} +\frac{-t^3(t^2 - t + 2)}{(t-1)(t^2+1)}$ 
\end{tabular}
\normalsize
\end{figure}

Examples for this theorem are given for the prime $\ell = 3$ and thirty-six Bianchi groups in figure~\ref{3torsionsubgraphs} (for $\ell = 2$, see \cite{RahmThesis}). In all the non-Euclidean principal ideal domain cases, the 2--torsion, and respectively 3--torsion subcomplexes are homeomorphic, which explains the results in proposition~\ref{non-Euclidean results}. 
Underlying theorem~\ref{invarianceUnderReduction}, there is the following correspondence between the non-trivial cyclic subgroups of the vertex stabilisers and the geodesic lines around which they effect a rotation, and which we shall call \emph{rotation axes}.

\begin{theorem}[~\cite{RahmThesis}] \label{actionOnAxes}
For any vertex $v$ in hyperbolic space, the action of its stabiliser on the set of rotation axes passing through $v$, induced by the action of the \mbox{Bianchi} group, is equivalent to the conjugation action of this stabiliser on its non-trivial cyclic subgroups.
\end{theorem}

A case by case study \cite{RahmThesis} for all the six types of finite sub-groups in the Bianchi groups allows us to deduce from theorem~\ref{actionOnAxes} the following lemma, which is useful in order to obtain theorem~\ref{invarianceUnderReduction}.

\begin{lemma} \label{geometricRigiditytheorem}
Let $v$ be a vertex in the refined cell complex. Then the number~$\bf n$ of orbits of edges in the refined cell complex adjacent to $v$, with stabiliser in ${\rm PSL}_2(\ringO_{-m})$ isomorphic to~$\Z/ \ell\Z$,  is given as follows for $\ell = 2$ and $\ell = 3$.
$$ \begin{array}{c|cccccc}
{\rm Isomorphism}\medspace { \rm type}\medspace {\rm of } \medspace {\rm the  }\medspace {\rm vertex }\medspace {\rm stabiliser} & \{1\} & \Z/2 \Z & \Z/3 \Z & \Z/2 \Z \times \Z/2 \Z & \Sthree & \Afour \\ 
\hline &&&&&& \scriptsize \\
{\bf n} \medspace \mathrm{ for } \medspace \ell = 2 & 0 & 2 & 0 & 3 & 2 & 1 \\ 
{\bf n} \medspace \mathrm{ for } \medspace \ell = 3 & 0 & 0 & 2 & 0 & 1 & 2.
\end{array} $$ \normalsize
\end{lemma}

Here we have written $\Sthree$ for the symmetric group on three letters and $\Afour$ for the alternating group on four letters. Note that we obtain the same table for ${\rm SL}_2(\ringO_{-m})$ after replacing the edge and vertex stabiliser types by their pre-images under the projection ${\rm SL}_2(\ringO_{-m}) \to {\rm PSL}_2(\ringO_{-m})$, which are respectively:  $\Z/2 \Z$, $\Z/4\Z$, $\Z/6\Z$, the 8-elements quaternion group, the 12-elements binary dihedral group and the binary tetrahedral group.

\subsection*{\emph{K}-theory} With the above information about the action of the \mbox{Bianchi} groups, we can further compute the Bredon homology of the \mbox{Bianchi} groups, from which we can deduce their equivariant $K$-homology. 
The results of the computations~\cite{RahmThesis} are the following.

\begin{theorem} \label{K-homology}
Let $\beta_1$ be the Betti number specified in proposition {\rm \ref{non-Euclidean results}}. For $\ringO_{-m}$ principal, the equivariant $K$-homology of $\Gamma := \mathrm{PSL}_2(\ringO_{-m})$ is isomorphic to \small
$$ \begin{array}{l|cccccc}

                                &m=1   & m=2            & m=3              & m=7  & m= 11            & m \in \{19,43,67,163\} \\ 

\hline &&&&&& \\

K^\Gamma_0(\underbar{E} \Gamma) & \Z^6 &\Z^5 \oplus \Z/2 \Z& \Z^5 \oplus \Z/2 \Z & \Z^3 & \Z^4 \oplus \Z/2 \Z & \Z^{\beta_1 -1} \oplus \Z^3 \oplus \Z/2 \Z \\ 

\\

K^\Gamma_1(\underbar{E} \Gamma) & \Z   & \Z^3           &  0               & \Z^3 &  \Z^3            & \Z \oplus \Z^{\beta_1}. \\ 

   \end{array}$$ \normalsize 
\end{theorem}

The remainder of the equivariant $K$-homology of $\Gamma$ is given by Bott $2$-periodicity.
By the Baum/Connes conjecture, which holds for the \mbox{Bianchi} groups \cite{JulgKasparov},
we obtain the $K$-theory of the reduced $C^*$-algebras of the \mbox{Bianchi} groups as isomorphic images.

\subsection*{Complex orbifolds} 
The information we have concerning the action of the \mbox{Bianchi} groups on real hyperbolic three-space provides explicit  orbifold structures.
These orbifolds serve as models in Cosmology \cite{AurichSteinerThen}.
\\
We can complexify these orbifolds. In \cite{RahmThesis}, the product structure on their Chen/Ruan orbifold cohomology has been determined for all \mbox{Bianchi} groups; and an algorithm has been given to compute the underlying vector space structure. This is a step towards checking Ruan's cohomological crepant resolution conjecture~\cite{ChenRuan}, which is of importance in Mathematical Physics and is completely open in 3 complex dimensions outside the global quotient case, hence these orbifolds are interesting test cases. 

\scriptsize
\textbf{Acknowledgements}. The author thanks Stephen Gelbart for support and encouragement. \\ He thanks \mbox{Philippe Elbaz-Vincent,} Louis Funar and Anthony Joseph for a careful lecture of the manuscript.
\\ And he would like to thank again all the people acknowledged in \cite{RahmThesis}.
\normalsize

\bibliographystyle{amsalpha}
\bibliography{literature}

\end{document}